\documentclass[a4paper,12pt]{article} 

\usepackage{color}

\usepackage{amsfonts, amsmath, amsthm, amssymb}
\usepackage[T1]{fontenc}
\usepackage[cp1250]{inputenc}
\usepackage{xcolor}
\usepackage{graphicx}
\usepackage{amssymb}
\usepackage{amsmath}
\usepackage{mathptmx}
\usepackage{helvet}
\usepackage{courier}
\usepackage{txfonts}
\usepackage{tikz} 
\usetikzlibrary{arrows}
\usepackage{type1cm}

\usepackage{verbatim}

\usepackage{graphicx}
\usepackage{epsfig,amscd,amssymb,amsxtra,amsmath,amsthm}
\usepackage{type1cm}
\usepackage[T1]{fontenc}
\usepackage{graphics}
\usepackage[mathscr]{eucal}
\usepackage[all]{xy}
\usepackage{amsmath,amscd}

\newtheorem{theorem}{Theorem}[section]

\newtheorem{definition}[theorem]{Definition}
\newtheorem{lemma}[theorem]{Lemma}
\newtheorem{claim}[theorem]{Claim}
\newtheorem{example}[theorem]{Example}
\newtheorem{remark}[theorem]{Remark}

\newtheorem{problem}[theorem]{Problem}
\newtheorem{observation}[theorem]{Observation}

\newcommand{\Cl}  {\mathop{\rm Cl}\nolimits}
\newcommand{\Int}  {\mathop{\rm Int}\nolimits}

\begin{document}

\def\joinrel{\mkern-3mu}
\newcommand{\varproj}{\displaystyle \lim_{\multimapinv\joinrel-\joinrel-}}

\title{The Lelek fan as the inverse limit of intervals with a single set-valued bonding function whose graph is an arc}
\author{Iztok Bani\v c,  Goran Erceg, and Judy Kennedy}
\date{}

\maketitle

\begin{abstract}
\noindent {  We consider a family of inverse limits of inverse sequences of closed unit intervals  with a single upper semi-continuous set-valued bonding function whose graph is an arc; it is  the union of two line segments in $[0,1]^2$, both of them contain the origin $(0, 0)$, have positive slope,  and extend to the opposite boundary of $[0,1]^2$. 

We show that there is a large subfamily $\mathcal F$ of these  bonding functions such that for each $f\in \mathcal F$, the inverse limit of the inverse sequence of closed unit intervals using $f$ as a single bonding function, is homeomorphic to the Lelek fan. }  
\end{abstract}
\-

\section{Introduction}

 {  Many famous continua such as the $\sin\frac{1}{x}$-continuum, the Knaster bucket-handle continuum, the pseudo-arc, and many more have been obtained as inverse limits of inverse sequences of closed unit intervals using a continuous single-valued function $[0,1]\rightarrow [0,1]$ as the only bonding function (see \cite{inmah} and \cite{nadler} for more such examples). Note that only chainable continua can be obtained in such a way, see \cite[Theorem 12.19, page 246]{nadler} for more details. Later, when a generalization of such inverse limits was obtained by Mahavier in \cite{mahavier} by introducing so-called generalized inverse limits or  inverse limits of inverse sequences of closed unit intervals using upper semi-continuous set-valued function from $[0,1]$ to the set of non-empty closed subsets of $[0,1]$, $[0,1]\multimap [0,1]$,  as the bonding functions, many new examples appeared. For example, the Wa\.{z}ewski's universal dendrite was obtained by Bani\v c, \v Crepnjak, Merhar, Milutinovi\' c and Sovi\v c  in \cite{banic} as the inverse limit of an inverse sequence of closed unit intervals using an upper semi-continuous set-valued function $[0,1]\multimap [0,1]$ as the only bonding function. Also,  the universal dendrite $D_m$, the Gehman dendrite, and the so-called monster continuum were obtained as such inverse limits in \cite{banic1,crepnjak} (by Bani\v c and Mart\' inez-de-la-Vega, and by \v Crepnjak and Sovi\v c), \cite{lemez} (by Leme\v z) and \cite[Example 2.15, page 31]{ingram} (by Ingram), respectively. More beautiful examples of continua that are presented in such an elegant way  may be found in \cite{ingram}.  However, not every continuum can be presented as such an inverse limit of an inverse sequence of closed unit intervals using an upper semi-continuous set-valued function $[0,1]\multimap [0,1]$ as the only bonding function.  It was proved in \cite{illanes} by Illanes that circle cannot be obtained in such a way, and in \cite{nall} it was proved by Nall that arcs are the only graphs that can be obtained as such inverse limits.  We conclude this list of various continua that can or cannot be obtained as inverse limits  of inverse sequences of closed unit intervals using an upper semi-continuous set-valued function $[0,1]\multimap [0,1]$ as the only bonding function by  a beautiful result from \cite{greenwood} that is obtained  by Greenwood and Suabedissen, showing that the only $2$-manifold that can be obtained as the inverse limit of an inverse sequence of closed unit intervals with upper semi-continuous set-valued bonding functions $[0,1]\multimap [0,1]$ is a $2$-torus. Then a generalization of this result was proved in \cite{greenwood1} by Alvin, Greenwood, Kelly and Kennedy that  the only $m$-manifold that can be obtained as such an inverse limit is an $m$-torus. 

In this paper, we construct a non-trivial family of fans as  the inverse limits of inverse sequences of closed unit intervals using an upper semi-continuous set-valued function $[0,1]\multimap [0,1]$ whose graph is an arc as the only bonding function. In particular, our main result is Theorem \ref{Lelek}, where the famous Lelek fan is presented as such an inverse limit. To obtain this surprising result, we consider a family of generalized inverse limits generated by two line segments in $[0,1]^2$. Both line segments contain the origin $(0,0)$, have positive slope, and extend to the opposite boundary of $[0,1]^2$. The generalized inverse limit is generated by the union of the two line segments. 

We proceed as follows. In Section \ref{s1}, the basic definitions and results that are needed later in the paper, are presented. In Section \ref{ss}, we introduce our settings and give a non-trivial family of fans as  the inverse limits of inverse sequences of closed unit intervals using an upper semi-continuous set-valued function $[0,1]\multimap [0,1]$ as the only bonding function.  One of the obtained fans is homeomorphic to the  Cantor fan. In Section \ref{s2}, our main results are presented; i.e., in Theorem \ref{Lelek},  we present the Lelek fan as  the inverse limit of an inverse sequence of closed unit intervals using an upper semi-continuous set-valued function $[0,1]\multimap [0,1]$ whose graph is an arc as the only bonding function.}

\section{Definitions and notation}\label{s1}
Suppose $X$ is a compact metric space. Recall that if $f:X \to X$ is a continuous function, the \emph{\color{blue} inverse limit space} generated by $f$ is 
\begin{equation*}
 \varprojlim(X,f):=\Big\{(x_{1},x_{2},x_{3},\ldots ) \in \prod_{i=1}^{\infty} X \ | \ 
\text{ for each positive integer } i,x_{i}= f(x_{i+1})\Big\},
\end{equation*}
\noindent which we can abbreviate as $\underset{\longleftarrow }{\lim }f$. The map $f$ on $X$ induces a natural homeomorphism $\sigma$ on $\underset{\longleftarrow }{\lim }f$, which is called the \emph{\color{blue}  shift map}, and is defined by 
$$\sigma(x_{1},x_{2},x_3, x_4,\ldots)=(x_{2},x_{3},x_4,\dots)$$ for each $(x_{1},x_{2},x_3, x_4,\ldots)$ in  $\underset{\longleftarrow }{\lim }f$.

Generalized inverse limits, or inverse limits with set-valued functions, are a generalization of (standard) inverse limits.  Here, rather than beginning with a continuous function $f$ from a compact metric space $X$ to itself, we begin with an upper semicontinuous function $f$ from $X$ to the non-empty closed subsets of $X$. The \emph{\color{blue} generalized inverse limit}, or the \emph{\color{blue} inverse limit with set-valued mappings}, associated with this set-valued function $f$ is the set 
\begin{equation*}
\varproj(X,f):=\Big\{(x_1,x_{2},x_3,\ldots )\in \prod_{i=1}^{\infty} X \ | \  
\text{ for each positive integer }i,x_{i}\in f(x_{i+1})\Big\},
\end{equation*}
\noindent which is a closed subspace of   $\Pi _{i =1}^{\infty}X$  endowed with the product topology.  Here again, the shift map $\sigma$ defined above takes $\varproj(X,f)$ to itself, but it is no longer a homeomorphism: $\sigma: \varproj(X,f) \to \varproj(X,f)$ is a continuous function.  And again, we often abbreviate $\varproj(X,f)$ as $\varproj f$. The topic of generalized inverse limits  is currently an
intensely studied area of continuum theory, with papers from many authors at this point. See \cite{ingram} for more references and additional information on the topic.
\begin{definition}
Let $(X,d)$ be a compact metric space. Then we define \emph{\color{blue}$2^X$} by 
$$
2^{X}=\{A\subseteq X \ | \ A \textup{ is a non-empty closed subset of } X\}.
$$
Let $\varepsilon >0$ and let $A\in 2^X$. Then we define  \emph{\color{blue}$N_d(\varepsilon,A)$} by 
$$
N_d(\varepsilon,A)=\bigcup_{a\in A}B(a,\varepsilon).
$$
Let $A,B\in 2^X$. The function \emph{\color{blue}$H_d:2^X\times 2^X\rightarrow \mathbb R$}, defined by
$$
H_d(A,B)=\inf\{\varepsilon>0 \ | \ A\subseteq N_d(\varepsilon,B), B\subseteq N_d(\varepsilon,A)\},
$$
is called \emph{\color{blue}a Hausdorff metric}. The Hausdorff metric is in fact a metric, the metric space $(2^X,H_d)$ is called \emph{\color{blue}the hyperspace of the space $(X,d)$}. 
\end{definition}
\begin{remark}
Let $(X,d)$ be a compact metric space, let $A$ be a non-empty closed subset of $X$,  and let $(A_n)$ be a sequence of non-empty closed subsets of $X$. When we say $\displaystyle A=\lim_{n\to \infty}A_n$ with respect to the Fausdorff metric, we mean $\displaystyle A=\lim_{n\to \infty}A_n$ in $(2^X,H_d)$. 
\end{remark}
\begin{definition}
Let $X$ and $Y$ be compact metric spaces. 
A function $F: X\rightarrow 2^Y$ is called \emph{ \color{blue}  a set-valued function} from $X$ to $Y$. We denote set-valued functions $F: X\rightarrow 2^Y$ by \emph{ \color{blue}  $F: X\multimap Y$}.
\end{definition}
\begin{definition}
A set-valued function  $F : X\multimap Y$ is  \emph{ \color{blue}   upper semicontinuous  at a point $x_0\in X$},  if for each
open set  $U\subseteq Y$ such that $F(x_0)\subseteq U$, there is an open set $V$ in $X$ such that
\begin{enumerate}
\item $x_0\in V$ and
\item for each $x\in V$, $F(x)\subseteq U$. 
\end{enumerate}  
The set-valued function  $F : X\multimap Y$ is  \emph{ \color{blue}   upper semicontinuous},  if it is upper semicontinuous at any point $x\in X$.
\end{definition}
\begin{definition}
 \emph{ \color{blue} The  graph $\Gamma(F)$ of a set-valued function} $F:X\multimap Y$ is the set of
all points  $(x,y)\in X\times Y$ such that $y \in F(x)$.  The set-valued function $F$ is surjective, if $\bigcup_{x\in X}F(x)=Y$. 
\end{definition}
There is a simple characterization of upper semicontinuous set-valued functions
(\cite[Proposition 11, p.\ 128]{A} and \cite[Theorem 1.2, p.\ 3]{ingram}):

\begin{theorem}
\label{th:grafi}  Let $X$ and $Y$ be compact metric spaces and $F:X\multimap Y$ a set-valued function. Then $F$ is upper semicontinuous if and only if its
graph $\Gamma(F)$ is closed in  $X\times Y$. 
\end{theorem}
\begin{definition}
Let $X$ be a compact metric space and let $G\subseteq X\times X$ be a relation on $X$. If $G$ is closed in $X\times X$, then we say that $G$ is  \emph{ \color{blue} a closed relation on $X$}.  
\end{definition}

\begin{definition}
Let $X$  be a set and let $G$ be a relation on $X$.  Then we define  
$$
G^{-1}=\{(y,x)\in X\times X \ | \ (x,y)\in G\}
$$
to be  \emph{ \color{blue} the inverse relation of the relation $G$ on $X$}.
\end{definition}
\begin{definition}
Let $X$ be a compact metric space and let $G$ be a closed relation on $X$. Then we call
$$
\star_{i=1}^{m}G^{-1}=\Big\{(x_1,x_2,x_3,\ldots ,x_{m+1})\in \prod_{i=1}^{m+1}X \ | \ \textup{ for each } i\in \{1,2,3,\ldots ,m\}, (x_{i+1},x_i)\in G\Big\}
$$
for each positive integer $m$,  \emph{ \color{blue} the $m$-th Mahavier product of $G$}, and
$$
\star_{i=1}^{\infty}G^{-1}=\Big\{(x_1,x_2,x_3,\ldots )\in \prod_{i=1}^{\infty}X \ | \ \textup{ for each positive integer } i, (x_{i+1},x_i)\in G\Big\}
$$
 \emph{ \color{blue} the infinite  Mahavier product of $G$}.
\end{definition}
\begin{observation}\label{obsi}
Let $X$ be a compact metric space, let $f:X\rightarrow X$ be a continuous function. 
Then 
$$
\star_{n=1}^{\infty}\Gamma(f)^{-1}=\varprojlim(X,f).
$$
Also, if $F:X\multimap X$ is an upper semi-continuous function, then  
$$
\star_{n=1}^{\infty}\Gamma(F)^{-1}=\varproj(X,f).
$$
\end{observation}

\begin{definition}
 \emph{ \color{blue} A continuum} is a non-empty compact connected metric space.  \emph{ \color{blue} A subcontinuum} is a subspace of a continuum, which is itself a continuum.
 \end{definition}
\begin{definition}
Let $X$ be a continuum. 
\begin{enumerate}
\item The continuum $X$ is \emph{\color{blue} unicoherent}, if for any subcontinua $A$ and $B$ of $X$ such that $X=A\cup B$,  the compactum $A\cap B$ is connected. 
\item The continuum $X$ is \emph{\color{blue} hereditarily unicoherent } provided that each of its subcontinua is unicoherent.
\item The continuum $X$ is a \emph{\color{blue} dendroid}, if it is an arcwise connected, hereditarily unicoherent continuum.
\item Let $X$ be a continuum.  If $X$ is homeomorphic to $[0,1]$, then $X$ is \emph{\color{blue} an arc}.   
\item A point $x$ in an arc $X$ is called \emph{\color{blue} an end-point of the arc  $X$}, if  there is a homeomorphism $\varphi:[0,1]\rightarrow X$ such that $\varphi(0)=x$.
\item Let $X$ be a dendroid.  A point $x\in X$ is called an \emph{\color{blue} end-point of the dendroid $X$}, if for  every arc $A$ in $X$ that contains $x$, $x$ is an end-point of $A$.  The set of all end-points of $X$ will be denoted by $E(X)$. 
\item A continuum $X$ is \emph{\color{blue} a simple triod}, if it is homeomorphic to $([-1,1]\times {0})\cup (\{0\}\times [0,1])$.
\item A point $x$ in a simple triod $X$ is called \emph{\color{blue} the top-point} or, breafly, the \emph{\color{blue} top of the simple triod $X$}, if  there is a homeomorphism $\varphi:([-1,1]\times {0})\cup (\{0\}\times [0,1])\rightarrow X$ such that $\varphi(0,0)=x$.
\item Let $X$ be a dendroid.  A point $x\in X$ is called \emph{\color{blue} a ramification-point of the dendroid $X$}, if there is a simple triod $T$ in $X$ with the top   $x$.  The set of all ramification-points of $X$ will be denoted by $R(X)$. 
\item The continuum $X$ is \emph{\color{blue} a  fan}, if it is a dendroid with at most one ramification point $v$, which is called the top of the fan $X$ (if it exists).
\item Let $X$ be a fan.   For all points $x$ and $y$ in $X$, we define  \emph{\color{blue} $A[x,y]$} to be the arc in $X$ with end-points $x$ and $y$, if $x\neq y$. If $x=y$, then we define $A[x,y]=\{x\}$.
\item Let $X$ be a fan with the top $v$. We say that that the fan $X$ is \emph{\color{blue} smooth} if for any $x\in X$ and for any sequence $(x_n)$ of points in $X$,
$$
\lim_{n\to \infty}x_n=x \Longrightarrow \lim_{n\to \infty}A[v,x_n]=A[v,x].
$$ 
\item Let $X$ be a fan.  We say that $X$ is \emph{\color{blue} a Cantor fan}, if $X$ is homeomorphic to the continuum
$$
\bigcup_{c\in C}A_c,
$$
where $C\subseteq [0,1]$ is the standard Cantor set and for each $c\in C$, $A_c$ is the straight line segment in the plane from $(0,0)$ to $(c,1)$.
\item Let $X$ be a fan.  We say that $X$ is \emph{\color{blue} a Lelek fan}, if it is smooth and $\Cl(E(X))=X$.
\end{enumerate}
\end{definition}
\begin{observation}
Let $X$ be a fan and let $Y$ be a subcontinuum of $X$. Then also $Y$ is a fan. 
\end{observation}
\begin{observation}
Suppose that  a fan $X$ with the top $v$ in the Hilbert cube $Q=\prod_{i=1}^{\infty}[0,1]$ is the union 
$$
X=\bigcup_{x\in E(X)}L_{x},
$$
where each $L_{x}$ is is a straight line segment (a convex segment) in $Q$ from $v$ to $x$. Then $X$ is a smooth fan. 
\end{observation}
\begin{observation}
Note that it was proved in \cite{lelek} by A.~Lelek that a Lelek fan exists.  Also, note that it has been proved in \cite{charatonik}  by W.~Charatonik and later,  in \cite{oversteegen} by W.~D.~Bula and L.~Overseegen  that arbitrary Lelek fans are homeomorphic. 
\end{observation}
\section{A  family of beautiful fans}\label{ss}
{  In this section, we give in Example \ref{beauty} a non-trivial family of fans presented as  inverse limits of inverse sequences of closed unit intervals using an upper semi-continuous set-valued function $[0,1]\multimap [0,1]$ as the only bonding function. Before we do that, we introduce the setting and the notation that will be used to present our results.} 
\begin{definition}
For each $(r,\rho)\in (0,\infty)\times (0,\infty)$, we define the sets \emph{\color{blue} $L_r$}, \emph{\color{blue} $L_{\rho}$} and \emph{\color{blue} $L_{r,\rho}$}  as follows:
$$
L_r=\{(x,y)\in [0,1]\times [0,1] \ | \ y=rx\},
$$
$$
L_{\rho}=\{(x,y)\in [0,1]\times [0,1] \ | \ y=\rho x\}
$$
and
$$
L_{r,\rho}=L_r\cup L_{\rho}.
$$
\end{definition}

\begin{definition}
For each $(r,\rho)\in (0,\infty)\times (0,\infty)$, we define the set \emph{\color{blue} $M_{r,\rho}$}  as follows:
$$
M_{r,\rho}=\star_{i=1}^{\infty}L_{r,\rho}.
$$
\end{definition}
\begin{definition}
Let $A$ be a set. We use  \emph{\color{blue} $A^{\mathbb N}$} to denote the set 
$$
\{a \ | \ a:\mathbb N\rightarrow A\}
$$
of all sequences in $A$.
\end{definition}
\begin{theorem}
For each $(r,\rho)\in (0,\infty)\times (0,\infty)$, 
$$
M_{r,\rho}=\bigcup_{a\in \{r,\rho\}^{\mathbb N}}\star_{i=1}^{\infty}L_{a_i}.
$$
\end{theorem}
\begin{proof}
Let $\mathbf x=(x_1,x_2,x_3,\ldots)\in M_{r,\rho}$. Then $(x_i,x_{i+1})\in L_{r,\rho}$ for each positive integer $i$. For each positive integer $i$, we let $a_i=r$, if $(x_i,x_{i+1})\in L_{r}$, and  $a_i=\rho$, if $(x_i,x_{i+1})\in L_{\rho}$. It follows that $\mathbf x\in \star_{i=1}^{\infty}L_{a_i}$. Next, let $\mathbf x\in \star_{i=1}^{\infty}L_{a_i}$ for some $a\in \{r,\rho\}^{\mathbb N}$. Since for each positive integer $i$, $L_{a_i}\subseteq L_{r,\rho}$, it follows that $\mathbf x\in M_{r,\rho}$.
\end{proof}
\begin{definition}
For each $(r,\rho)\in (0,\infty)\times (0,\infty)$ and for each sequence $a\in \{r,\rho\}^{\mathbb N}$, we define the set \emph{\color{blue} $L_{a}$}  as follows:
$$
L_{a}=\star_{i=1}^{\infty}L_{a_i}.
$$
For each $(r,\rho)\in (0,\infty)\times (0,\infty)$, we also define $\mathcal L_{r,\rho}=\{L_a \ | \ a\in \{r,\rho\}^{\mathbb N}\}$.
\end{definition}
\begin{example}\label{beauty}
Let $(r,\rho)\in (0,\infty)\times (0,\infty)$. To describe $M_{r,\rho}$,  we consider the following possible cases.
\begin{enumerate}
\item $r=\rho$. In this case, $M_{r,\rho}=\star_{i=1}^{\infty}L_{r}$. In addition, 
\begin{enumerate}
\item if $r\leq 1$, then $M_{r,\rho}$ is an arc with end points $(0,0,0,\ldots )$ and $(1,r,r^2,r^3,\ldots)$.
\item if $r>1$, then $M_{r,\rho}=\{(0,0,0,\ldots )\}$.
\end{enumerate}
\item $r\neq \rho$.  We consider the following subcases.
\begin{enumerate}
\item $r>1$ and $\rho >1$. Here, $M_{r,\rho}=\{(0,0,0,\ldots )\}$.
\item $r\leq 1$ and $\rho \leq 1$. For each sequence  $a\in \{r,\rho\}^{\mathbb N}$, $L_a$ is a straight line segment in the Hilbert cube $Q=\prod_{i=1}^{\infty}[0,1]$ with endpoints $(0,0,0,\ldots)$ and $(1,a_1,a_2\cdot a_1, a_3\cdot a_2\cdot a_1,\ldots)$.  Note that $\{r,\rho\}^{\mathbb N}$ is a Cantor set and that 
$$
\varphi:\{r,\rho\}^{\mathbb N}\rightarrow \{(1,a_1,a_2\cdot a_1, a_3\cdot a_2\cdot a_1,\ldots) \ | \ a\in \{r,\rho\}^{\mathbb N}\},
$$
defined by 
$$
\varphi(a)=(1,a_1,a_2\cdot a_1, a_3\cdot a_2\cdot a_1,\ldots)
$$
for any $a\in \{r,\rho\}^{\mathbb N}$, is a homeomorphism.  Therefore, 
$$
\{(1,a_1,a_2\cdot a_1, a_3\cdot a_2\cdot a_1,\ldots) \ | \ a\in \{r,\rho\}^{\mathbb N}\}
$$
is a Cantor set in $Q$.  Since $M_{r,\rho}=\bigcup_{a\in \{r,\rho\}^{\mathbb N}}L_a$, it follows that $M_{r,\rho}$ is a Cantor fan. 
\item $r=1$ and $\rho > 1$. Note that for each sequence $a\in \{r,\rho\}^{\mathbb N}$, 
$$
L_a=\{(0,0,0,\ldots)\} \Longleftrightarrow a_i=\rho \textup{ for infinitely many indexes } i.
$$
For each non-negative integer $n$, let 
$$
A_n=\{a\in \{r,\rho\}^{\mathbb N} \ | \ a_i=\rho \textup{ for exactly } n \textup{ indexes } i\},
$$ 
let $\mathcal M_n=\{L_a \ | \ a\in A_n\}$ and let $M_n=\bigcup_{a\in A_n}L_a$.
Then $M_{r,\rho}=\bigcup_{n=0}^{\infty}M_n$. Note that $\mathcal M_0=\{L_{(1,1,1,\ldots)}\}$ and that for each positive integer $n$, $\mathcal M_n$ is a collection of countable many arcs, each of them having diameter less or equal to 
$$
D_n=\max\left\{\frac{1}{\rho^n}, \frac{1}{2\cdot \rho^{n-1}}, \frac{1}{2^2\cdot \rho^{n-2}}, \ldots, \frac{1}{2^{n-1}\cdot \rho},\frac{1}{2^n}\right\}
$$ 
and each of them having $(0,0,0,\ldots)$ as one of their endpoints.
Also, note that $\displaystyle \lim_{n\to \infty}D_n=0$ and that for each non-negative integer $n$, for each $a\in A_n$, and for each positive integer $m > n$, there is a sequence $(L_{i})$ of arcs in $\mathcal M_{m}$ such that for for some subarc $L$ of $L_a$, $\displaystyle \lim_{i\to\infty}L_i=L$ with respect to the Hausdorff metric.  
In particular, $M_0$ is the strait line segment in the Hilbert cube $Q$ with end points $(0,0,0,\ldots)$ and $(1,1,1,\ldots )$, while the set $M_1=\bigcup_{i=1}^{\infty}L_i$,  where for each positive integer $i$,
$$
L_i=\left\{(\underbrace{t,t,t,\ldots ,t}_{i},\rho\cdot t, \rho\cdot t,\rho\cdot t,\ldots) \ | \ t\in \Big[0,\frac{1}{\rho}\Big]\right\};
$$
see Figure \ref{figure1}. 
 \begin{figure}[h!]
	\centering
		\includegraphics[width=35em]{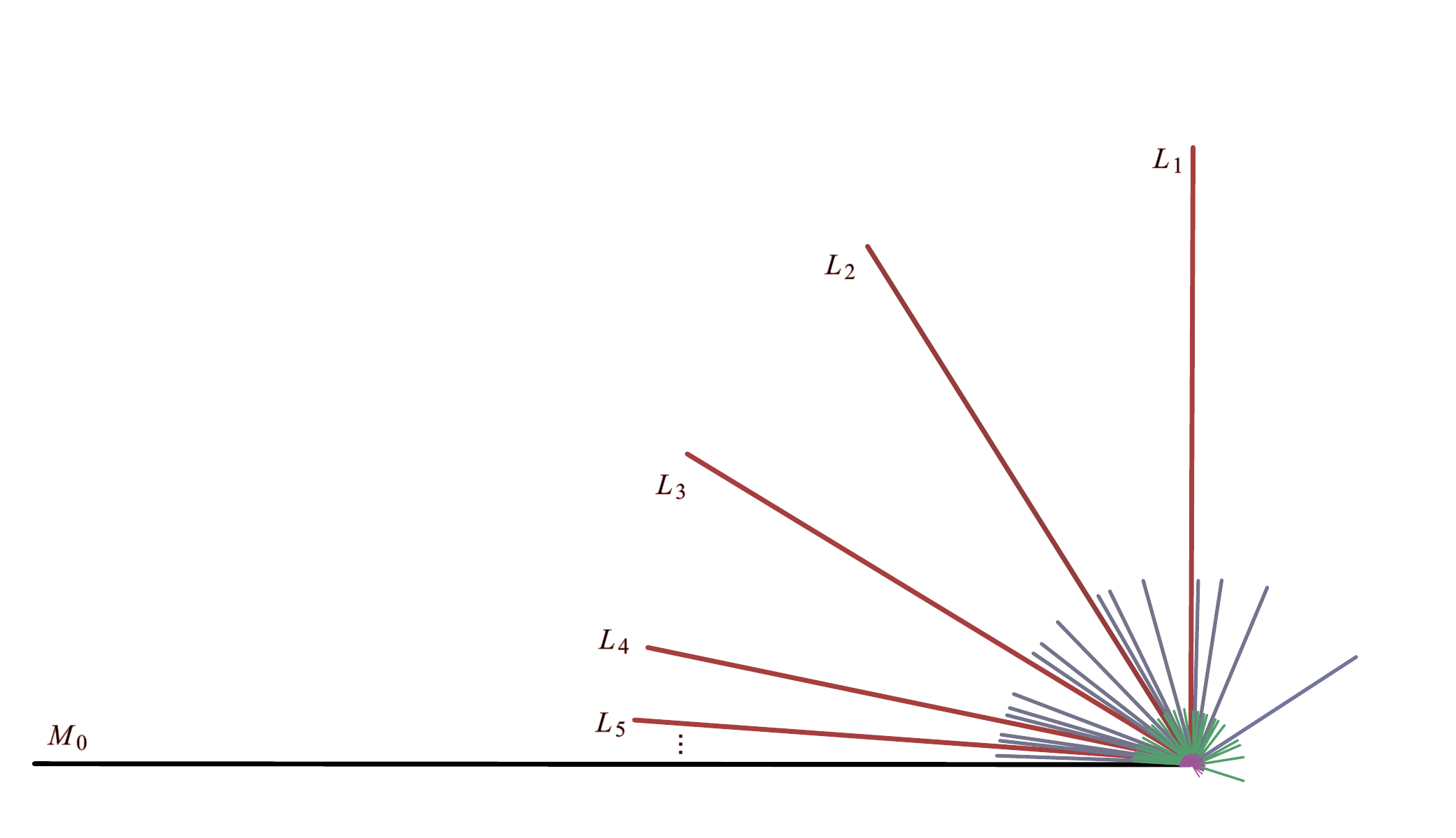}
	\caption{The fan $M_{r,\rho}$ for $\rho=1$ and $r > 1$}
	\label{figure1}
\end{figure} 
\item $r<1$, $\rho>1$ and for all integers $k$ and $\ell$,  
$$
r^k = \rho^{\ell} \Longleftrightarrow k=\ell=0.
$$
Section \ref{s2} is dedicated to showing in Theorem \ref{Lelek} that in this case,  $M_{r,\rho}$ is homeomorphic to the Lelek fan. 
\item $r<1$, $\rho>1$ and there are integers $k$ and $\ell$ such that 
\begin{enumerate}
\item $k\neq 0$ or $\ell\neq 0$, and
\item  $r^k = \rho^{\ell}$.
\end{enumerate}
We leave this case as an open problem; see Problem \ref{problem}.
\end{enumerate}
\end{enumerate} 
\end{example}
\section{The Lelek fan}\label{s2}
{  In this section, our main results are presented. Before stating and proving them, we give the following definitions.}
\begin{definition}
Let $r,\rho\in \mathbb R$. We say that \emph{\color{blue} $r$ and $\rho$ never connect} or \emph{\color{blue} $(r,\rho)\in \mathcal{NC}$}, if \begin{enumerate}
\item $r<1$, $\rho>1$ and 
\item for all integers $k$ and $\ell$,  
$$
r^k = \rho^{\ell} \Longleftrightarrow k=\ell=0.
$$
\end{enumerate} 
\end{definition}

\begin{definition}
Let  $a,b\in \mathbb R$ be such that $a<b$ and let $t\in \mathbb R\setminus \{0\}$. Then we define 
$$
t\cdot (a,b)=(ta,tb).
$$
\end{definition}
{  In the theorems that follow, we prove the fundamental properties of $\mathcal{NC}$ that will be used to prove Theorem \ref{Lelek}, which is our main result.}
\begin{theorem}\label{dense1}
Let  $(r,\rho)\in \mathcal{NC}$.  The set
$$
\{r^k\cdot \rho^{\ell} \ | \ k \textup{ and } \ell \textup{ are  integers}\}
$$
is dense in $(0,\infty)$.
\end{theorem}
\begin{proof}
Let 
$$
\mathcal B=\{r^k\cdot \rho^{\ell} \ | \ k \textup{ and } \ell \textup{ are  integers}\}
$$
and suppose that $\mathcal B$ is not dense in $(0,\infty)$.  Then there are $\alpha,\beta \in [0,\infty)$ such that
\begin{enumerate}
\item $\alpha < \beta$,
\item $(\alpha,\beta)\cap \mathcal B=\emptyset$,  and
\item \label{tri} for all $\gamma, \delta\in [0,\infty)$, 
$$
(\alpha,\beta)\subseteq (\gamma,\delta) \textup{ and }  (\gamma,\delta)\cap \mathcal B=\emptyset \Longrightarrow \alpha=\gamma \textup{ and } \beta=\delta.
$$
\end{enumerate}
Choose such numbers $\alpha$ and $\beta$.  Note that it follows that $\alpha,\beta\in \Cl(\mathcal B)$. 
\begin{claim}\label{1ena}
For all integers $k$ and $\ell$, 
$$
r^k\rho^{\ell}\cdot (\alpha,\beta)\cap \mathcal B=\emptyset. 
$$
\end{claim}
\begin{proof}[Proof of Claim \ref{1ena}]\renewcommand{\qed}{} 
Suppose that for integers $k_0,\ell_0$, $r^{k_0}\rho^{\ell_0}\in r^k\rho^{\ell}\cdot (\alpha,\beta)$.  It follows that $r^{k_0-k}\rho^{\ell_0-\ell}\in (\alpha,\beta)$---a contradiction. This completes the proof of the claim. 
\end{proof}
\begin{claim}\label{2dva}
For all integers $k$, $\ell$, $k'$ and $\ell'$, 
$$
r^k\rho^{\ell}\neq r^{k'}\rho^{\ell'}  \Longrightarrow  (r^k\rho^{\ell}\cdot (\alpha,\beta))\cap (r^{k'}\rho^{\ell'}\cdot (\alpha,\beta))=\emptyset. 
$$
\end{claim}
\begin{proof}[Proof of Claim \ref{2dva}]\renewcommand{\qed}{} 
Let $k$, $\ell$, $k'$ and $\ell'$ be such that $r^k\rho^{\ell}\neq r^{k'}\rho^{\ell'}$ and suppose that $(r^k\rho^{\ell}\cdot (\alpha,\beta))\cap (r^{k'}\rho^{\ell'}\cdot (\alpha,\beta))\neq\emptyset$. Let 
$$
z\in (r^k\rho^{\ell}\cdot (\alpha,\beta))\cap (r^{k'}\rho^{\ell'}\cdot (\alpha,\beta)).
$$
It follows that 
$$
r^{-k}\rho^{-\ell}\cdot z\in (\alpha,\beta)\cap (r^{k'-k}\rho^{\ell'-\ell}\cdot (\alpha,\beta)).
$$
Since $r^{-k}\rho^{-\ell}\cdot z\in (\alpha,\beta)$ and since $(\alpha,\beta)\cap \mathcal B=\emptyset$, it follows that $r^{-k}\rho^{-\ell}\cdot z\not \in \mathcal B$.  Also, since $(\alpha,\beta)$ and $(r^{k'-k}\rho^{\ell'-\ell}\cdot (\alpha,\beta))$ are both open intervals, each of them containing the point $r^{-k}\rho^{-\ell}\cdot z$, it follows that $(\alpha,\beta)\cup (r^{k'-k}\rho^{\ell'-\ell}\cdot (\alpha,\beta))$ is also an open interval in $\mathbb R$.   It follows from \ref{tri} that  $r^{k'-k}\rho^{\ell'-\ell}\cdot (\alpha,\beta)\subseteq (\alpha,\beta)$.  Therefore,
$$
(\alpha,\beta)\subseteq r^{k-k'}\rho^{\ell-\ell'}\cdot (\alpha,\beta).
$$
By Claim \ref{1ena},  $r^{k-k'}\rho^{\ell-\ell'}\cdot (\alpha,\beta)\cap \mathcal B=\emptyset$. and it follows from \ref{tri} that 
$$
(\alpha,\beta)=r^{k-k'}\rho^{\ell-\ell'}\cdot (\alpha,\beta).
$$
Therefore,  $r^{k-k'}\rho^{\ell-\ell'}=1$ and 
$$
r^{k'-k}=\rho^{\ell-\ell'}
$$
follows.  Hence, $k'=k$ and $\ell'=\ell$, meaning that $r^k\rho^{\ell} = r^{k'}\rho^{\ell'}$---a contradiction.
  This completes the proof of the claim. 
\end{proof}
Next, we define a sequence $z:\mathbb N\rightarrow [0,1]\cap \mathcal B$ by defining $z_1=1$, $z_2=r$, and for each positive integer $n\geq 2$,
$$
z_{n+1}=\begin{cases}
				\rho\cdot z_n\text{;} & \rho\cdot z_n\leq 1 \\
				r\cdot z_n\text{;} & \rho\cdot z_n>1.
			\end{cases}
$$
\begin{claim}\label{3tri}
For each positive integer $n$, 
$$
\frac{r}{\rho}\leq z_n\leq 1.
$$
\end{claim}
\begin{proof}[Proof of Claim \ref{1ena}]\renewcommand{\qed}{} 
It follows from the construction of the sequence $z$ that for each positive integer $n$, 
$z_n\leq 1$.  It follows from $r<\rho$ that $1\geq \frac{r}{\rho}$ and $r\geq \frac{r}{\rho}$.  Therefore, $z_1\geq \frac{r}{\rho}$ and $z_2\geq \frac{r}{\rho}$. Let $n$ be a positive integer and suppose that $z_n\geq \frac{r}{\rho}$. To prove that $z_{n+1}\geq  \frac{r}{\rho}$,  we consider the following cases.
\begin{enumerate}
\item $z_{n+1}=\rho\cdot z_n$. Then 
$$
z_{n+1}=\rho\cdot z_n\geq z_n\geq  \frac{r}{\rho}
$$
and we are done.
\item $z_{n+1}=r\cdot z_n$.  Then $\rho\cdot z_n>1$ and it follows that $z_n>\frac{1}{\rho}$. Therefore,
$$
z_{n+1}=r\cdot z_n>r\cdot \frac{1}{\rho}= \frac{r}{\rho}.
$$
\end{enumerate}
This completes the proof of the claim. 
\end{proof}
Next, let 
$$
\mathcal A=\{z_n\cdot (\alpha,\beta) \ | \  n \textup{ is a positive integer}\}.
$$
It follows from Claim \ref{2dva} that for all positive integers $m$ and $n$,
$$
m\neq n \Longrightarrow (z_m\cdot (\alpha,\beta))\cap (z_n\cdot (\alpha,\beta))=\emptyset.
$$
Also, note that for each positive integer $n$, 
\begin{enumerate}
\item $z_n\cdot (\alpha,\beta)\subseteq [0,\beta]$, since $0<z_n\leq 1$, and 
\item $z_n\beta-z_n\alpha=z_n(\beta-\alpha)\geq \frac{r}{\rho}\cdot (\beta-\alpha)$.
\end{enumerate} 
Therefore, there are infinitely many pairwise disjoint open intervals of length at least  $\frac{r}{\rho}\cdot (\beta-\alpha)$ in the interval $[0,\beta]$, which is not possible. This completes the proof.
\end{proof}
\begin{definition}
Let $(r,\rho)\in \mathcal{NC}$. We define 
$$
\mathcal B_1(r,\rho)=\{r^k\cdot \rho^{\ell} \ | \ k \textup{ and } \ell \textup{ are  non-negative integers}\},
$$
$$
\mathcal B_2(r,\rho)=\{r^k\cdot \rho^{\ell} \ | \ k \textup{ and } \ell \textup{ are  negative integers}\},
$$
$$
\mathcal B_3(r,\rho)=\{r^k\cdot \rho^{\ell} \ | \ k \textup{ is a  non-negative integer} \textup{ and } \ell \textup{ is a negative integer}\},
$$
$$
\mathcal B_4(r,\rho)=\{r^k\cdot \rho^{\ell} \ | \ k \textup{ is a  negative integer} \textup{ and } \ell \textup{ is a non-negative integer}\},
$$
and
$$
\mathcal B(r,\rho)=\mathcal B_1(r,\rho)\cup \mathcal B_2(r,\rho)\cup\mathcal B_3(r,\rho)\cup \mathcal B_4(r,\rho).
$$
\end{definition}
\begin{lemma}\label{nowhere}
Let $(r,\rho)\in \mathcal{NC}$.  Then $\mathcal B_3(r,\rho)$ and $\mathcal B_4(r,\rho)$ are nowhere dense in $(0,\infty)$. 
\end{lemma}
\begin{proof}
Note that for each $x\in (0,\infty)$,
\begin{enumerate}
\item $\mathcal B_3(r,\rho)\cap (x,\infty)$ is finite, and that
\item $\mathcal B_4(r,\rho)\cap (0,x)$ is finite. 
\end{enumerate}
Therefore,  both, $\mathcal B_3(r,\rho)$ and $\mathcal B_4(r,\rho)$, are closed in $(0,\infty)$ and neither of them contains an open interval.  It follows that $\mathcal B_3(r,\rho)$ and $\mathcal B_4(r,\rho)$ are nowhere dense in $(0,\infty)$. 
\end{proof}
\begin{lemma}\label{between}
Let $(r,\rho)\in \mathcal{NC}$.  Then  for all $x,y\in \mathcal B_2(r,\rho)$ , 
$$
x<y \Longrightarrow (x,y)\cap \mathcal B_1(r,\rho)\neq \emptyset.
$$
\end{lemma}
\begin{proof}
Let $x,y\in \mathcal B_2(r,\rho)$ such that $x<y$.  Also, let $k_1$, $k_2$, $\ell_1$ and $\ell_2$ be positive integers such that 
$$
x=r^{-k_1}\rho^{-\ell_1} \textup{ and  } y=r^{-k_2}\rho^{-\ell_2}. 
$$
Note that 
\begin{enumerate}
\item $(0,x)\cap \mathcal B_1(r,\rho)\neq \emptyset$ and 
\item $(y,\infty)\cap \mathcal B_1(r,\rho)\neq \emptyset$. 
\end{enumerate}
We consider the following cases.
\begin{enumerate}
\item \label{ABBA} $k_1=k_2$.  First, note that 
$$
x<y\Longleftrightarrow r^{-k_1}\rho^{-\ell_1} <r^{-k_2}\rho^{-\ell_2} \Longleftrightarrow \rho^{-\ell_1} <\rho^{-\ell_2}\Longleftrightarrow -\ell_1<-\ell_2 \Longleftrightarrow \ell_2<\ell_1.
$$
Suppose that $ (x,y)\cap \mathcal B_1(r,\rho)=\emptyset$.  First, we prove the following claim. 
\begin{claim}\label{one}
There are a positive integer $k$ and a non-negative integer  $\ell$ such that  
$$
r^{k}\rho^{\ell}<x<y<r^{k}\rho^{\ell+1}.
$$
\end{claim}
\begin{proof}[Proof of Claim \ref{one}]\renewcommand{\qed}{} 
Let $k$ be a positive integer such that $r^k\leq x$.  Since $\rho>1$, it follows that \begin{enumerate}
\item for each non-negative integer $n$, $r^{k}\rho^{n}<r^{k}\rho^{n+1}$ and
\item $\displaystyle \lim_{n\to\infty} r^{k}\rho^{n}=\infty$.
\end{enumerate} 
Therefore,  there is a non-negative integer $\ell$ such that $r^{k}\rho^{\ell}\leq x$ and $r^{k}\rho^{\ell+1}> x$.  Since $ (x,y)\cap \mathcal B_1(r,\rho)=\emptyset$, it follows that $r^{k}\rho^{\ell+1}\geq y$.  All that is left to see is that $r^{k}\rho^{\ell}\neq x$ and $r^{k}\rho^{\ell+1}\neq y$.   Suppose that $r^{k}\rho^{\ell}= x$. It follows that $r^{k}\rho^{\ell}= r^{-k_1}\rho^{-\ell_1}$.  Therefore, $r^{k+k_1}=\rho^{-\ell_1-\ell}$ and this implies that $k+k_1=0$ and $-\ell_1-\ell=0$. Therefore, $k=-k_1$---a contradiction, since $k_1$ and $k$ are both positive integers. Suppose that $r^{k}\rho^{\ell+1}= y$. It follows that $r^{k}\rho^{\ell+1}= r^{-k_2}\rho^{-\ell_2}$.  Therefore, $r^{k+k_2}=\rho^{-\ell_2-\ell-1}$ and this implies that $k+k_2=0$ and $-\ell_2-\ell-1=0$. Therefore, $k=-k_2$---a contradiction, since $k_2$ and $k$ are both positive integers. This completes the proof of the claim. 
\end{proof}
Let $k$ be a positive integer  and let $\ell$ be a non-negative integer such that  
$$
r^{k}\rho^{\ell}<x<y<r^{k}\rho^{\ell+1}.
$$
Then 
$$
r^{-k_1}\rho^{-\ell_2}=r^{-k_2}\rho^{-\ell_2}=y<r^{k}\rho^{\ell+1}=r^{k}\rho^{\ell}\cdot \rho<x\cdot \rho= r^{-k_1}\rho^{-\ell_1}\cdot \rho=r^{-k_1}\rho^{-\ell_1+1}.
$$
It follows that $r^{-k_1}\rho^{-\ell_2}<r^{-k_1}\rho^{-\ell_1+1}$. 
Note that
$$
r^{-k_1}\rho^{-\ell_2}<r^{-k_1}\rho^{-\ell_1+1} \Longleftrightarrow \rho^{-\ell_2}<\rho^{-\ell_1+1} \Longleftrightarrow -\ell_2<-\ell_1+1 \Longleftrightarrow \ell_1-1<\ell_2. 
$$
It follows that $\ell_1\leq \ell_2$---a contradiction, since $\ell_2<\ell_1$. Therefore,  $ (x,y)\cap \mathcal B_1(r,\rho)\neq \emptyset$.  
\item $\ell_1=\ell_2$.  Note that
$$
x<y\Longleftrightarrow r^{-k_1}\rho^{-\ell_1} <r^{-k_2}\rho^{-\ell_2} \Longleftrightarrow r^{-k_1} <r^{-k_2}\Longleftrightarrow -k_1>-k_2 \Longleftrightarrow k_1<k_2.
$$
Suppose that $ (x,y)\cap \mathcal B_1(r,\rho)=\emptyset$.  First, we prove the following claim. 
\begin{claim}\label{two}
There are a non-negative integer $k$ and a positive integer  $\ell$ such that  
$$
r^{k+1}\rho^{\ell}<x<y<r^{k}\rho^{\ell}.
$$
\end{claim}
\begin{proof}[Proof of Claim \ref{two}]\renewcommand{\qed}{} 
Let $\ell$ be a positive integer such that $\rho^{\ell}\geq y$.  Since $r<1$, it follows that \begin{enumerate}
\item for each non-negative integer $n$, $r^{n+1}\rho^{\ell}<r^{n}\rho^{\ell}$ and
\item $\displaystyle \lim_{n\to\infty} r^{n}\rho^{\ell}=0$.
\end{enumerate} 
Therefore,  there is a non-negative integer $k$ such that $r^{k}\rho^{\ell}\geq y$ and $r^{k+1}\rho^{\ell}< y$.  Since $ (x,y)\cap \mathcal B_1(r,\rho)=\emptyset$, it follows that $r^{k+1}\rho^{\ell}\leq x$.  All that is left to see is that $r^{k+1}\rho^{\ell}\neq x$ and $r^{k}\rho^{\ell}\neq y$.   Suppose that $r^{k+1}\rho^{\ell}= x$. It follows that $r^{k+1}\rho^{\ell}= r^{-k_1}\rho^{-\ell_1}$.  Therefore, $r^{k+1+k_1}=\rho^{-\ell_1-\ell}$ and this implies that $k+1+k_1=0$ and $-\ell_1-\ell=0$. Therefore, $\ell=-\ell_1$---a contradiction, since $\ell_1$ and $\ell$ are both positive integers.  Next, suppose that $r^{k}\rho^{\ell}= y$. It follows that $r^{k}\rho^{\ell}= r^{-k_2}\rho^{-\ell_2}$.  Therefore, $r^{k+k_2}=\rho^{-\ell_2-\ell}$ and this implies that $k+k_2=0$ and $-\ell_2-\ell=0$. Therefore, $\ell=-\ell_2$---a contradiction, since $\ell_2$ and $\ell$ are both positive integers. This completes the proof of the claim. 
\end{proof}
Let $k$ be a non-negative integer  and let $\ell$ be a positive integer such that  
$$
r^{k+1}\rho^{\ell}<x<y<r^{k}\rho^{\ell}.
$$
Then 
$$
r^{-k_2+1}\rho^{-\ell_1}=r\cdot r^{-k_2}\rho^{-\ell_1}=r\cdot r^{-k_2}\rho^{-\ell_2}=r\cdot y<r\cdot r^{k}\rho^{\ell}=r^{k+1}\rho^{\ell}<x= r^{-k_1}\rho^{-\ell_1}.
$$
It follows that $r^{-k_2+1}\rho^{-\ell_1}<r^{-k_1}\rho^{-\ell_1}$. 
Note that
$$
r^{-k_2+1}\rho^{-\ell_1}<r^{-k_1}\rho^{-\ell_1} \Longleftrightarrow r^{-k_2+1}<r^{-k_1} \Longleftrightarrow -k_2+1>-k_1 \Longleftrightarrow k_2-1<k_1. 
$$
It follows that $k_2\leq k_1$---a contradiction, since $k_1<k_2$. Therefore,  $ (x,y)\cap \mathcal B_1(r,\rho)\neq\emptyset$ follows. 
\item $k_1<k_2$. First, we prove the following claim. 
\begin{claim}\label{three}
If $\ell_2<\ell_1$, then $ (x,y)\cap \mathcal B_1(r,\rho)\neq\emptyset$. 
\end{claim}
\begin{proof}[Proof of Claim \ref{three}]\renewcommand{\qed}{} 
Suppose that  $\ell_2<\ell_1$. It follows that $\rho^{\ell_1}>\rho^{\ell_2}$ and it follows from this that $\rho^{-\ell_1}<\rho^{-\ell_2}$.  Note that $r^{k_2-k_1}<1$. Therefore,
 $$
 y= r^{-k_2}\rho^{-\ell_2}>r^{-k_2}\rho^{-\ell_2}\cdot r^{k_2-k_1}=r^{-k_1}\rho^{-\ell_2}>r^{-k_1}\rho^{-\ell_1}=x.
 $$
 Let $y'=r^{-k_1}\rho^{-\ell_2}$ and let $k'=k_1$ and $\ell'=\ell_2$.   We have 
 $$
 x=r^{-k_1}\rho^{-\ell_1} \textup{ and } y'=r^{-k'}\rho^{-\ell'}.
 $$
 Note that $x<y'$ and that $k_1=k'$. In the first case, we have seen that in this case,  $(x,y')\cap\mathcal B_1(r,\rho)\neq\emptyset $. Therefore $(x,y)\cap \mathcal B_1(r,\rho)\neq\emptyset$. This completes the proof of the claim. 
\end{proof}
Suppose until the end of this case that $\ell_1<\ell_2$ (note that we have already covered the cases $\ell_2<\ell_1$ and $\ell_1=\ell_2$).
\begin{claim}\label{fouri}
$r^{k_1-k_2}\rho^{\ell_1-\ell_2}>1$.
\end{claim}
\begin{proof}[Proof of Claim \ref{fouri}]\renewcommand{\qed}{} 
The claim follows from $x<y$ and 
$$
x\cdot r^{k_1-k_2}\rho^{\ell_1-\ell_2}=r^{-k_1}\rho^{-\ell_1}\cdot r^{k_1-k_2}\rho^{\ell_1-\ell_2}=r^{-k_2}\rho^{-\ell_2}=y.
$$
This completes the proof of the claim.
\end{proof}
\begin{claim}\label{four}
For each  positive integer $\ell$, 
$$
\rho^{\ell}\not \in \{y\cdot r^{n(k_1-k_2)}\rho^{n(\ell_1-\ell_2)} \ | \ n \textup{ is a positive integer}\}.
$$
\end{claim}
\begin{proof}[Proof of Claim \ref{four}]\renewcommand{\qed}{} 
Suppose that $\ell$ and $n$ are positive integers such that 
$$
\rho^{\ell}=y\cdot r^{n(k_1-k_2)}\rho^{n(\ell_1-\ell_2)}. 
$$
It follows that
$$
\rho^{\ell}=r^{-k_2}\rho^{-\ell_2}\cdot r^{n(k_1-k_2)}\rho^{n(\ell_1-\ell_2)}. 
$$
Therefore,
$$
r^{n(k_1-k_2)-k_2}=\rho^{\ell+\ell_2-n(\ell_1-\ell_2)}.
$$
Since $(r,\rho)\in \mathcal{NC}$, it follows that $n(k_1-k_2)-k_2=0$ and $\ell+\ell_2-n(\ell_1-\ell_2)=0$. Therefore, $k_2=n(k_1-k_2)$. Since $k_1-k_2<0$, it follows that $k_2<0$---a contradiction.  This completes the proof of the claim. 
\end{proof}
Let $\ell$ be a positive integer such that $\rho^{\ell}>y$.  By Claim \ref{four}, 
$$
\rho^{\ell}\in (y,\infty)\setminus \{y\cdot r^{n(k_1-k_2)}\rho^{n(\ell_1-\ell_2)} \ | \ n \textup{ is a positive integer}\}.
$$
It follows from Claim \ref{fouri} that for each non-negative integer $n$, 
$$
y\cdot r^{n(k_1-k_2)}\rho^{n(\ell_1-\ell_2)}\geq y,
$$
$$
y\cdot r^{(n+1)(k_1-k_2)}\rho^{(n+1)(\ell_1-\ell_2)} > y\cdot r^{n(k_1-k_2)}\rho^{n(\ell_1-\ell_2)},
$$
and
$$
\displaystyle \lim_{n\to\infty}y\cdot r^{n(k_1-k_2)}\rho^{n(\ell_1-\ell_2)}=\infty.
$$
Therefore,  there is a non-negative integer $n$ such that 
$$
\rho^{\ell}\in (y\cdot r^{n(k_1-k_2)}\rho^{n(\ell_1-\ell_2)} ,y\cdot r^{(n+1)(k_1-k_2)}\rho^{(n+1)(\ell_1-\ell_2)}).
$$
Let $n$ be such an integer and let 
$$
z= \rho^{\ell}\cdot r^{(n+1)(k_2-k_1)}\rho^{(n+1)(\ell_2-\ell_1)}.
$$ 
To complete the proof of this case, we prove the following claim.
\begin{claim}\label{five}
$z\in (x,y)\cap \mathcal B_1$.
\end{claim}
\begin{proof}[Proof of Claim \ref{five}]\renewcommand{\qed}{} 
To show that $z\in \mathcal B_1$, observe that
$$
z= \rho^{\ell}\cdot r^{(n+1)(k_2-k_1)}\rho^{(n+1)(\ell_2-\ell_1)}=r^{(n+1)(k_2-k_1)}\rho^{(n+1)(\ell_2-\ell_1)+\ell}.
$$
Since $k_2>k_1$ and $\ell_2>\ell_1$, it follows that $(n+1)(k_2-k_1)\geq 0$ and $(n+1)(\ell_2-\ell_1)+\ell\geq 0$.  Therefore, $z\in \mathcal B_1$.  To see that $z\in (x,y)$, observe that it follows from 
$$
\rho^{\ell}\in (y\cdot r^{n(k_1-k_2)}\rho^{n(\ell_1-\ell_2)} ,y\cdot r^{(n+1)(k_1-k_2)}\rho^{(n+1)(\ell_1-\ell_2)} )
$$
 and 
 $$
 z= \rho^{\ell}\cdot r^{(n+1)(k_2-k_1)}\rho^{(n+1)(\ell_2-\ell_1)}
 $$
  that 
$$
{\color{blue} r^{(n+1)(k_2-k_1)}\rho^{(n+1)(\ell_2-\ell_1)}\cdot y\cdot r^{n(k_1-k_2)}\rho^{n(\ell_1-\ell_2)}}<z
$$
and 
$$
z<{\color{red} r^{(n+1)(k_2-k_1)}\rho^{(n+1)(\ell_2-\ell_1)}\cdot y\cdot r^{(n+1)(k_1-k_2)}\rho^{(n+1)(\ell_1-\ell_2)})}.
$$
Since 
$$
{\color{blue} r^{(n+1)(k_2-k_1)}\rho^{(n+1)(\ell_2-\ell_1)}\cdot y\cdot r^{n(k_1-k_2)}\rho^{n(\ell_1-\ell_2)}}=r^{-k_1}\rho^{-\ell_1}=x
$$
and 
$$
{\color{red} r^{(n+1)(k_2-k_1)}\rho^{(n+1)(\ell_2-\ell_1)}\cdot y\cdot r^{(n+1)(k_1-k_2)}\rho^{(n+1)(\ell_1-\ell_2)})}=y,
$$
it follows that $z\in (x,y)$.  This completes the proof of the claim. 
\end{proof}
\item $k_1>k_2$.  First, we prove the following claim.  
\begin{claim}\label{six}
$\ell_2<\ell_1$.
\end{claim}
\begin{proof}[Proof of Claim \ref{six}]\renewcommand{\qed}{} 
It follows from $k_1>k_2$ that $r^{-k_1}>r^{-k_2}$.  Therefore, $0<r^{k_1-k_2}<1$. Observe that
$$
x<y \Longleftrightarrow r^{-k_1}\rho^{-\ell_1}<r^{-k_2}\rho^{-\ell_2} \Longleftrightarrow \rho^{\ell_2-\ell_1}<r^{k_1-k_2}.
$$
Therefore, $\rho^{\ell_2-\ell_1}<1$. It follows that $\ell_2<\ell_1$. This completes the proof of the claim.
\end{proof}
\begin{claim}\label{fourii}
$r^{k_2-k_1}\rho^{\ell_2-\ell_1}<1$.
\end{claim}
\begin{proof}[Proof of Claim \ref{fourii}]\renewcommand{\qed}{} 
The claim follows from $x<y$ and 
$$
y\cdot r^{k_2-k_1}\rho^{\ell_2-\ell_1}=r^{-k_2}\rho^{-\ell_2}\cdot r^{k_2-k_1}\rho^{\ell_2-\ell_1}=r^{-k_1}\rho^{-\ell_1}=x.
$$
This completes the proof of the claim.
\end{proof}
\begin{claim}\label{seven}
For each positive integer $k$, 
$$
r^k\not \in \{x\cdot r^{n(k_2-k_1)}\rho^{n(\ell_2-\ell_1)} \ | \ n \textup{ is a positive integer}\}.
$$
\end{claim}
\begin{proof}[Proof of Claim \ref{seven}]\renewcommand{\qed}{} 
Suppose that $k$ and $n$ are positive integers such that 
$$
r^k=x\cdot r^{n(k_2-k_1)}\rho^{n(\ell_2-\ell_1)}.
$$
It follows that
$$
r^k=r^{-k_1}\rho^{-\ell_1}\cdot r^{n(k_2-k_1)}\rho^{n(\ell_2-\ell_1)}.
$$
Therefore,
$$
r^{k+k_1-n(k_2-k_1)}=\rho^{n(\ell_2-\ell_1)-\ell_1}.
$$
Since $(r,\rho)\in \mathcal{NC}$, it follows that
$k+k_1-n(k_2-k_1)=0$ and $n(\ell_2-\ell_1)-\ell_1=0$. Therefore, $\ell_1=n(\ell_2-\ell_1)$. Since $\ell_2-\ell_1<0$, it follows that $\ell_1<0$---a contradiction. This completes the proof of the claim.
\end{proof}
Let $k$ be a positive integer such that $r^k<x$.  By Claim \ref{seven}, 
$$
r^k\not \in \{x\cdot r^{n(k_2-k_1)}\rho^{n(\ell_2-\ell_1)} \ | \ n \textup{ is a positive integer}\}.
$$
It follows from Claim \ref{fourii} that for each non-negative integer $n$, 
$$
x\cdot r^{n(k_2-k_1)}\rho^{n(\ell_2-\ell_1)} \leq x,
$$
$$
x\cdot r^{(n+1)(k_2-k_1)}\rho^{(n+1)(\ell_2-\ell_1)} < x\cdot r^{n(k_2-k_1)}\rho^{n(\ell_2-\ell_1)},
$$
and
$$
\displaystyle \lim_{n\to\infty}x\cdot r^{n(k_2-k_1)}\rho^{n(\ell_2-\ell_1)}=0.
$$
Therefore,  there is a non-negative integer $n$ such that 
$$
r^{k}\in (x\cdot r^{(n+1)(k_2-k_1)}\rho^{(n+1)(\ell_2-\ell_1)} ,x\cdot r^{n(k_2-k_1)}\rho^{n(\ell_2-\ell_1)}).
$$
Let $n$ be such an integer and let 
$$
z= r^k\cdot r^{(n+1)(k_1-k_2)}\rho^{(n+1)(\ell_1-\ell_2)}.
$$ 
To complete the proof of this case, we prove the following claim.
\begin{claim}\label{nine}
$z\in (x,y)\cap \mathcal B_1$.
\end{claim}
\begin{proof}[Proof of Claim \ref{nine}]\renewcommand{\qed}{} 
To show that $z\in \mathcal B_1$, observe that
$$
z=r^k\cdot r^{(n+1)(k_1-k_2)}\rho^{(n+1)(\ell_1-\ell_2)}=r^{k+(n+1)(k_1-k_2)}\rho^{(n+1)(\ell_1-\ell_2)}.
$$
Since $k_2<k_1$ and $\ell_2<\ell_1$, it follows that $k+(n+1)(k_1-k_2)\geq 0$ and $(n+1)(\ell_1-\ell_2)\geq 0$.  Therefore, $z\in \mathcal B_1$.  To see that $z\in (x,y)$, observe that it follows from 
$$
r^{k}\in (x\cdot r^{(n+1)(k_2-k_1)}\rho^{(n+1)(\ell_2-\ell_1)} ,x\cdot r^{n(k_2-k_1)}\rho^{n(\ell_2-\ell_1)})
$$
 and 
 $$
 z= r^k\cdot r^{(n+1)(k_1-k_2)}\rho^{(n+1)(\ell_1-\ell_2)}
 $$
  that 
$$
{\color{blue} r^{(n+1)(k_1-k_2)}\rho^{(n+1)(\ell_1-\ell_2)}\cdot x\cdot r^{(n+1)(k_2-k_1)}\rho^{(n+1)(\ell_2-\ell_1)}}<z
$$
and 
$$
z<{\color{red} r^{(n+1)(k_1-k_2)}\rho^{(n+1)(\ell_1-\ell_2)}\cdot x\cdot r^{n(k_2-k_1)}\rho^{n(\ell_2-\ell_1)})}.
$$
Since 
$$
{\color{blue} r^{(n+1)(k_1-k_2)}\rho^{(n+1)(\ell_1-\ell_2)}\cdot x\cdot r^{(n+1)(k_2-k_1)}\rho^{(n+1)(\ell_2-\ell_1)}}=x
$$
and 
$$
{\color{red} r^{(n+1)(k_1-k_2)}\rho^{(n+1)(\ell_1-\ell_2)}\cdot x\cdot r^{n(k_2-k_1)}\rho^{n(\ell_2-\ell_1)})}=r^{-k_2}\rho^{-\ell_2}=y,
$$
it follows that $z\in (x,y)$.  
\end{proof}
\end{enumerate}
This completes the proof. 
\end{proof}
In the proof of Theorem \ref{enkica}, we also use  the following well-known lemmas.  Since the proofs are short, we give them here for the completeness of the paper.
\begin{lemma}\label{iztok}
Let $X$ be a  metric space and $A,B\subseteq X$ such that $A$ is dense in $X$
 and $B$ is nowhere dense in $X$. Then $A\setminus B$ is dense in $X$.
 \end{lemma}
\begin{proof}
Let $U$ be a non-empty open set in $X$.  Since $B$ is nowhere dense in $X$, it follows that $\Int(\Cl(B))=\emptyset$. Therefore $U\not \subseteq \Cl(B)$ and it follows that $U\setminus \Cl(B)\neq \emptyset$.  Since $U\setminus \Cl(B)$ is non-empty and open in $X$ and since $A$ is dense in $X$, it follows that $A\cap (U\setminus \Cl(B))\neq \emptyset$.  Therefore,  $U\cap (A\setminus B)\neq \emptyset$.
\end{proof}
\begin{lemma}\label{iktoz}
Let $A,B\subseteq (0,\infty)$ be such sets that $A\cup B$ is dense in $X$, $A\cap B=\emptyset$  and for all $x,y\in B$,
 $$
 x<y \Longrightarrow (x,y)\cap A\neq \emptyset.
 $$ 
  Then $A$ is dense in $X$.
\end{lemma}
\begin{proof}
Let $a,b\in \mathbb R$ be such that $a<b$. Since $A\cup B$ is dense in $(0,\infty)$, it follows that $(a,b)\cap (A\cup B)\neq \emptyset$. Let $x\in (a,b)\cap (A\cup B)$.  Since $A\cup B$ is dense in $(0,\infty)$, it follows that $(x,b)\cap (A\cup B)\neq \emptyset$.  Let $y\in (x,b)\cap (A\cup B)$.  If $x\in A$ or $y\in A$, then we are done.  If not, then $x,y\in B$.  By the assumption,  $(x,y)\cap A\neq \emptyset$. Therefore, $(a,b)\cap A\neq \emptyset$. 
\end{proof}
\begin{lemma}\label{iktoz1}
Let $X$ be a  metric space and $A,B\subseteq X$ such that $A$ is nowhere dense in $X$
 and $B$ is nowhere dense in $X$. Then $A\cup B$ is nowhere dense in $X$.
\end{lemma}
\begin{proof}
To prove that $A\cup B$ is nowhere dense in $X$, we show that $X\setminus (\Cl(A\cup B))$ is dense in $X$.  Let $U$ be a non-empty open set in $X$.  Since $A$ and $B$ are nowhere dense in $X$, it follows that $X\setminus (\Cl(A))$ and $X\setminus (\Cl(B))$ are dense in $X$. Also, note that $X\setminus (\Cl(A))$ and $X\setminus (\Cl(B))$ are both open in $X$.  Therefore, since $U\cap (X\setminus (\Cl(A)))\neq \emptyset$ and since $X\setminus (\Cl(B))$  is dense in $X$, it follows that $U\cap (X\setminus (\Cl(A)))\cap (X\setminus (\Cl(B)))\neq \emptyset$. Therefore,
$$
\emptyset \neq U\cap (X\setminus (\Cl(A)))\cap (X\setminus (\Cl(B)))=U\cap (\Int(X\setminus A))\cap (\Int(X\setminus B))=
$$
$$
U\cap (\Int((X\setminus A)\cap (X\setminus B)))=U\cap \Int(X\setminus (A\cup B))=U\cap (X\setminus (\Cl(A\cup B)))
$$
It follows that $X\setminus (\Cl(A\cup B))$ is dense in $X$. 
\end{proof}
\begin{theorem}\label{enkica}
Let $(r,\rho)\in \mathcal{NC}$.  Then $\mathcal B_1(r,\rho)$ is dense in $(0,\infty)$. 
\end{theorem}
\begin{proof}
Recall that by Theorem \ref{dense1},
$$
\mathcal B(r,\rho)=\mathcal B_1(r,\rho)\cup \mathcal B_2(r,\rho)\cup\mathcal B_3(r,\rho)\cup \mathcal B_4(r,\rho).
$$
is dense in $(0,\infty)$. By Lemma \ref{nowhere}, $\mathcal B_3(r,\rho)$ and $\mathcal B_4(r,\rho)$ are nowhere dense in $(0,\infty)$ and by Lemma \ref{iktoz1}, $\mathcal B_3(r,\rho)\cup \mathcal B_4(r,\rho)$ is nowhere dense in $(0,\infty)$.  It follows from Lemma \ref{iztok} that 
$$
\mathcal B(r,\rho)\setminus (\mathcal B_3(r,\rho)\cup \mathcal B_4(r,\rho))=\mathcal B_1(r,\rho)\cup \mathcal B_2(r,\rho)
$$
is dense in $(0,\infty)$.  By Theorem \ref{between}, for all $x,y\in  \mathcal B_2(r,\rho)$,
 $$
 x<y \Longrightarrow (x,y)\cap  \mathcal B_1(r,\rho)\neq \emptyset.
 $$ 
 Since $\mathcal B_1(r,\rho)\cap \mathcal B_2(r,\rho)=\emptyset$, it follows from Lemma \ref{iktoz} that $\mathcal B_1(r,\rho)$ is dense in $(0,\infty)$. 
\end{proof}
The following lemma is also a very known result. Since the proof is short, we give it here for the completeness of the paper. 
\begin{lemma}\label{tale}
Let $A$ be dense in $(0,\infty)$ and let $x\in (0,\infty)$. Then $\{x\cdot a \ | \ a\in A\}$ is also dense in $(0,\infty)$. 
\end{lemma}
\begin{proof}
Let $\alpha,\beta\in (0,\infty)$ such that $\alpha<\beta$.  Since $A$ is dense in $(0,\infty)$, it follows that $A\cap (\frac{1}{x}\cdot \alpha,\frac{1}{x}\cdot \beta)\neq \emptyset$.  Let $a\in A\cap (\frac{1}{x}\cdot \alpha,\frac{1}{x}\cdot \beta)$.  Then $x\cdot a\in (\alpha, \beta)$. It follows that $\{x\cdot a \ | \ a\in A\}$ is also dense in $(0,\infty)$. 
\end{proof}
We use the following lemma in the proof of Theorem \ref{to}.
\begin{lemma}\label{goran}
Let $(r,\rho)\in \mathcal{NC}$. Then for each $x\in (0,1)$, there are sequences $(k_n)$ and $(\ell_n)$ of non-negative integers such that for each positive integer $n$, 
$$
1-\frac{1}{n+1}<x\cdot r^{k_n}\rho^{\ell_n}<1.
$$
\end{lemma}
\begin{proof}
By Theorem \ref{enkica},  $\mathcal B_1(r,\rho)$ is dense in $(0,\infty)$. By Lemma \ref{tale}, 
$$
\{x\cdot y \ | \ y\in \mathcal B_1(r,\rho)\}
$$
 is also dense in $(0,\infty)$. Therefore, for each positive integer $n$, there are non-negative integers $k_n$ and $\ell_n$ such that 
 $$
1-\frac{1}{n+1}<x\cdot r^{k_n}\rho^{\ell_n}<1.
 $$
\end{proof}
\begin{theorem}\label{to}
Let $(r,\rho)\in \mathcal{NC}$. Then for each $x\in (0,1)$,  there is a sequence $a\in \{r,\rho\}^{\mathbb N}$ such that for each positive integer $n$,
$$
(a_1\cdot a_2\cdot a_3\cdot \ldots \cdot a_n)\cdot x \in [0,1]
$$
and 
$$
\sup\{(a_1\cdot a_2\cdot a_3\cdot \ldots \cdot a_n)\cdot x \ | \ n \textup{ is a positive integer}\}=1.
$$
\end{theorem}
\begin{proof}
Let $x\in (0,1)$ and let $(k_n)$ and $(\ell_n)$ be sequences of non-negative integers such that for each positive integer $n$, 
$$
1-\frac{1}{n+1}<x\cdot r^{k_n}\rho^{\ell_n}<1.
$$
Such sequences do exist by Lemma \ref{goran}.  Then, let 
$$
a_1=a_2=a_3=\ldots =a_{k_1}=r \textup{ and } a_{k_1+1}=a_{k_1+2}=a_{k_1+3}=\ldots =a_{k_1+\ell_1}=\rho. 
$$
Since $\rho>1$ and $(a_1\cdot a_2\cdot a_3\cdot \ldots \cdot a_{k_1+\ell_1})\cdot x \in [0,1]$, it follows that for each $i\in 1,2,3,\ldots, k_1+\ell_1$,
$$
(a_1\cdot a_2\cdot a_3\cdot \ldots \cdot a_i)\cdot x \in [0,1].
$$
Next, we define 
$$
a_{k_1+\ell_1+1}=a_{k_1+\ell_1+2}=a_{k_1+\ell_1+3}=\ldots =a_{k_1+\ell_1+k_2}=r
$$
and
$$
 a_{k_1+\ell_1+k_2+1}=a_{k_1+\ell_1+k_2+2}=a_{k_1+\ell_1+k_2+3}=\ldots =a_{k_1+\ell_1+k_2+\ell_2}=\rho. 
$$
Obviously,  for each $i\in 1,2,3,\ldots, k_1+\ell_1+k_2+\ell_2$,
$$
(a_1\cdot a_2\cdot a_3\cdot \ldots \cdot a_i)\cdot x \in [0,1].
$$
Let $n$ be a positive integer and suppose that we have already defined the terms $a_1$, $a_2$, $a_3$, $\ldots$, $a_{\sum_{j=1}^{n}(k_j+\ell_j)}$ such that 
for each $i\in 1,2,3,\ldots, \sum_{j=1}^{n}(k_j+\ell_j)$,
$$
(a_1\cdot a_2\cdot a_3\cdot \ldots \cdot a_i)\cdot x \in [0,1].
$$ 
Then we define 
$$
a_{\sum_{j=1}^{n}(k_j+\ell_j)+1}=a_{\sum_{j=1}^{n}(k_j+\ell_j)+2}=a_{\sum_{j=1}^{n}(k_j+\ell_j)+3}=\ldots =a_{\sum_{j=1}^{n}(k_j+\ell_j)+k_{n+1}}=r
$$
and
$$
 a_{\sum_{j=1}^{n}(k_j+\ell_j)+k_{n+1}+1}=a_{\sum_{j=1}^{n}(k_j+\ell_j)+k_{n+1}+2}=a_{\sum_{j=1}^{n}(k_j+\ell_j)+k_{n+1}+3}=\ldots =
 $$
 $$
 a_{\sum_{j=1}^{n}(k_j+\ell_j)+k_{n+1}+\ell_{n+1}}=\rho. 
$$
Obviously,  for each $i\in 1,2,3,\ldots, \sum_{j=1}^{n+1}(k_j+\ell_j)$,
$$
(a_1\cdot a_2\cdot a_3\cdot \ldots \cdot a_i)\cdot x \in [0,1].
$$
We have just constructed (inductively) our sequence $a$.  Obviously, 
for each positive integer $n$,
$$
(a_1\cdot a_2\cdot a_3\cdot \ldots \cdot a_n)\cdot x \in [0,1].
$$
Note that
$$
\sup\{(a_1\cdot a_2\cdot a_3\cdot \ldots \cdot a_n)\cdot x \ | \ n \textup{ is a positive integer}\}=1
$$ 
follows from the fact that for each positive integer $n$, 
$$
1-\frac{1}{n+1}<x\cdot r^{k_n}\rho^{\ell_n}<1.
$$
This completes our proof.
\end{proof}

\begin{observation}
For each $(r,\rho)\in \mathcal{NC}$ and for each $a\in \{r,\rho\}^{\mathbb N}$, $L_a$ is either a straight line segment in the Hilbert cube $Q$ with one endpoint being the point  $(0,0,0,\ldots)$ or  $L_a=\{(0,0,0,\ldots)\}$. 
\end{observation}
\begin{definition}
Let  $(r,\rho)\in \mathcal{NC}$  and let $a\in  \{r,\rho\}^{\mathbb N}$. We say that the sequence $a$ is \emph{\color{blue} a useful sequence}, if $L_a$ is an arc.  We also define the set \emph{\color{blue} $\mathcal U_{r,\rho}$} as follows:
$$
\mathcal U_{r,\rho}=\{a\in  \{r,\rho\}^{\mathbb N} \ | \ a \textup{ is a useful sequence}\}.
$$
\end{definition}
\begin{definition}
 For each positive integer $k$,  we use $\pi_k:\prod_{i=1}^{\infty}[0,1]\rightarrow [0,1]$ to denote \emph{ \color{blue}  the $k$-th standard projection} from $\prod_{i=1}^{\infty}[0,1]$ to $[0,1]$.
\end{definition}
\begin{definition}
Let  $(r,\rho)\in \mathcal{NC}$. For each $a\in \mathcal U_{r,\rho}$, we define the sequence \emph{\color{blue} $b^a:\mathbb N\rightarrow [0,1]$} by 
$$
b^a_n=\max(\pi_n(L_a))
$$
for each positive integer $n$.
\end{definition}
\begin{observation}
Let  $(r,\rho)\in \mathcal{NC}$. For each $a\in \mathcal U_{r,\rho}$,  
$$
\pi_n(L_a)=[0,b^a_n]
$$
for each positive integer $n$. 
\end{observation}
{  Next, we show that $M_{r,\rho}$ is a fan.}
\begin{theorem}\label{fen}
Let  $(r,\rho)\in \mathcal{NC}$.  Then $M_{r,\rho}$ is a fan such that
$$
E(M_{r,\rho})=\{(b_1^a,b_2^a,b_3^a,\ldots) \ | \  a\in \mathcal U_{r,\rho}\} \textup{ and } R(M_{r,\rho})=\{(0,0,0,\ldots)\}.
$$ 
\end{theorem}
\begin{proof}
Let
$$
E=L_{r,\rho}\cup \Big(\Big[\frac{1}{\rho},1\Big]\times \{1\}\Big)
$$
and let 
$$
F=\star_{i=1}^{\infty}E.
$$
Note that $F$ is a Cantor fan:  it is the union  
$$
F=\bigcup_{a\in  \{r,1\}^{\mathbb N} }A_a,
$$
where for each $ a\in \{r,1\}^{\mathbb N}$,  $A_a$ is an arc with end-points $(0,0,0,\ldots)$ and $(a_1,a_1a_2,a_1a_2a_3,\ldots )$, also,  
$$
\{(a_1,a_1a_2,a_1a_2a_3,\ldots ) \ | \ a\in \{r,1\}^{\mathbb N} \}
$$
is homeomorphic to $\{r,1\}^{\mathbb N}$, which is a Cantor set.  

Since every subcontinuum of a fan is itself a fan, it follows from $M_{r,\rho}\subseteq F$,  that  $M_{r,\rho}$ is a fan.  Obviously, 
$$
E(M_{r,\rho})=\{(b_1^a,b_2^a,b_3^a,\ldots) \ | \  a\in \mathcal U_{r,\rho}\}
$$
and 
$$
R(M_{r,\rho})=\{(0,0,0,\ldots)\}.
$$
and this completes the proof. 
\end{proof}
{ In the following theorem, the end-points of the fan $M_{r,\rho}$ are found.}
\begin{theorem}\label{toto}
Let $(r,\rho)\in \mathcal{NC}$,  let $x\in (0,1)$ and let $a\in \{r,\rho\}^{\mathbb N}$ be a sequence such that for each positive integer $n$,
$$
(a_1\cdot a_2\cdot a_3\cdot \ldots \cdot a_n)\cdot x \in [0,1]
$$
and 
$$
\sup\{(a_1\cdot a_2\cdot a_3\cdot \ldots \cdot a_n)\cdot x \ | \ n \textup{ is a positive integer}\}=1.
$$  
Then 
$$
(x,a_1\cdot x, (a_1\cdot a_2)\cdot x,(a_1\cdot a_2\cdot a_3)\cdot x,\ldots )\in E(M_{r,\rho}).
$$
\end{theorem}
\begin{proof}
Suppose that $(x,a_1\cdot x, (a_1\cdot a_2)\cdot x,(a_1\cdot a_2\cdot a_3)\cdot x,\ldots )\not \in E(M_{r,\rho})$.  Let $a=(a_1,a_2,a_3,\ldots)$.  Then $x<b_1^a$.  It follows that for each positive integer $n$,
$$
(a_1\cdot a_2\cdot a_3\cdot \ldots \cdot a_n)\cdot x<(a_1\cdot a_2\cdot a_3\cdot \ldots \cdot a_n)\cdot b_1^a.
$$ 
Let $(i_n)$ be a strictly increasing sequence of positive integers such that 
$$
\lim_{n\to \infty}(a_1\cdot a_2\cdot a_3\cdot \ldots \cdot a_{i_n})\cdot x=1.
$$
Then 
$$
\lim_{n\to \infty}(a_1\cdot a_2\cdot a_3\cdot \ldots \cdot a_{i_n})\cdot b_1^a=1.
$$
It follows that
$$
0=1-1=\lim_{n\to \infty}(a_1\cdot a_2\cdot a_3\cdot \ldots \cdot a_{i_n})\cdot b_1^a-\lim_{n\to \infty}(a_1\cdot a_2\cdot a_3\cdot \ldots \cdot a_{i_n})\cdot x=
$$
$$
\lim_{n\to \infty}(a_1\cdot a_2\cdot a_3\cdot \ldots \cdot a_{i_n})\cdot (b_1^a-x)=(b_1^a-x)\lim_{n\to \infty}(a_1\cdot a_2\cdot a_3\cdot \ldots \cdot a_{i_n})
$$
and since $b_1^a-x\neq 0$, it follows that 
$$
\lim_{n\to \infty}(a_1\cdot a_2\cdot a_3\cdot \ldots \cdot a_{i_n})=0,
$$
which is a contradiction.  Therefore,  
$$
(x,a_1\cdot x, (a_1\cdot a_2)\cdot x,(a_1\cdot a_2\cdot a_3)\cdot x,\ldots )\in E(M_{r,\rho}).
$$ 
This completes the proof.
\end{proof}
The following theorem is our main result.
\begin{theorem}\label{Lelek}
Let $(r,\rho)\in \mathcal{NC}$. Then $M_{r,\rho}$ is a Lelek fan.
\end{theorem}
\begin{proof}
By Theorem \ref{fen},  $M_{r,\rho}$ is a fan.  Since it is the union 
$$
M_{r,\rho}=\bigcup_{a\in \mathcal{U}_{r,\rho}}L_a
$$
of straight line segments $L_a$ from $(0,0,0,\ldots)$ to $(b_1^a, b_2^a, b_3^a,\ldots )$ in the Hilbert cube $Q$, it is smooth.  To complete the proof, we prove that $E(M_{r,\rho})$ is dense in $M_{r,\rho}$.  It suffices to prove that $E(M_{r,\rho})$ is dense in $M_{r,\rho}\setminus \{(0,0,0,\ldots)\}$.  Let 
$$
\mathbf x=(x_1,x_2,x_3,\ldots)\in M_{r,\rho}\setminus \{(0,0,0,\ldots)\}
$$
be any point.  Then $x_n\neq 0$ for each positive integer $n$.  For each positive integer $n$, by Theorem \ref{to}, there is a sequence 
$$
a^n=(a_1^n, a_2^n, a_3^n, \ldots)\in \{r,\rho\}^{\mathbb N}
$$
 such that for each positive integer $k$,
$$
(a_1^n\cdot a_2^n\cdot a_3^n\cdot \ldots \cdot a_k^n)\cdot x_n \in [0,1]
$$
and 
$$
\sup\{(a_1^n\cdot a_2^n\cdot a_3^n\cdot \ldots \cdot a_k^n)\cdot x_n \ | \ k \textup{ is a positive integer}\}=1.
$$
For each positive integer $n$, choose such a sequence  $a^n$ and let
$$
\mathbf x_n=(x_1,x_2,x_3,\ldots , x_n,a_1^n\cdot x_n, (a_1^n\cdot a_2^n)\cdot x_n,(a_1^n\cdot a_2^n\cdot a_3^n)\cdot x_n,\ldots ).
$$
By Theorem \ref{toto},  $\mathbf x_n \in E(M_{r,\rho})$ for each positive integer $n$. Since 
$$
\lim_{n\to \infty}\mathbf x_n=\mathbf x,
$$
it follows that $E(M_{r,\rho})$ is dense in $M_{r,\rho}\setminus \{(0,0,0,\ldots)\}$.  Therefore,  $E(M_{r,\rho})$ is dense in $M_{r,\rho}$.   
\end{proof}
\begin{observation}
Let $(r, \rho) \in \mathcal{NC}$ and let $f : [0, 1] \multimap[0,1]$ be a set-valued function such that the graph of $f$ equals to  $\Gamma(f) = L_{r,\rho}^{-1}$. The set $L_{r,\rho}^{-1}$ is a closed subset of  $[0, 1] \times [0, 1]$, therefore, $f$ is an upper semi-continuous function. By Observation \ref{obsi} and by Theorem \ref{Lelek}, the generalized inverse limit $\varproj(X,f)$ is homeomorphic to the Lelek fan.
\end{observation}
{  We have just presented the Lelek fan as the inverse limit of an inverse sequence of closed unit intervals using an upper semi-continuous set-valued function $[0,1]\multimap [0,1]$ whose graph is an arc as the only bonding function.  Therefore, this is a good place where we can state an open problem and finish the paper.}
\begin{problem}\label{problem}
Let $r,r'<1$, $\rho,\rho'>1$ and let $k$,  $\ell$, $k'$ and $\ell'$ be integers  such that 
\begin{enumerate}
\item $k\neq 0$ or $\ell\neq 0$, and
\item  $r^k = \rho^{\ell}$,
\item $k'\neq 0$ or $\ell'\neq 0$, and
\item  ${(r')}^{k'} = {(\rho')}^{\ell'}$.
\end{enumerate}
Are the continua $M_{r,\rho}$ and $M_{r',\rho'}$ homeomorphic?
\end{problem}

\noindent I. Bani\v c\\
              (1) Faculty of Natural Sciences and Mathematics, University of Maribor, Koro\v{s}ka 160, SI-2000 Maribor,
   Slovenia; \\(2) Institute of Mathematics, Physics and Mechanics, Jadranska 19, SI-1000 Ljubljana, 
   Slovenia; \\(3) Andrej Maru\v si\v c Institute, University of Primorska, Muzejski trg 2, SI-6000 Koper,
   Slovenia\\
             {iztok.banic@um.si}           
     
				\-
				
		\noindent G.  Erceg\\
             Faculty of Science, University of Split, Rudera Bo\v skovi\' ca 33, Split,  Croatia\\
{{gorerc@pmfst.hr}       }    

                 \-

                 	\-
					
  \noindent J.  Kennedy\\
             Lamar University, 200 Lucas Building, P.O. Box 10047, Beaumont, TX 77710 USA\\
{{kennedy9905@gmail.com}       }    

\end{document}